\documentclass[11pt]{amsart}
\usepackage{amsmath}

\begin{document}

\newtheorem{thm}{Theorem}[section]
\newtheorem{lem}[thm]{Lemma}
\newtheorem{prop}[thm]{Proposition}
\newtheorem{cor}[thm]{Corollary}
\newtheorem{defn}[thm]{Definition}
\newtheorem*{remark}{Remark}

\numberwithin{equation}{section}

\newcommand{\Z}{{\mathbb Z}} %cph changed from \mathbf
\newcommand{\Q}{{\mathbb Q}}
\newcommand{\R}{{\mathbb R}}
\newcommand{\C}{{\mathbb C}}
\newcommand{\N}{{\mathbb N}}
\newcommand{\FF}{{\mathbb F}}

\newcommand{\rmk}[1]{\footnote{{\bf Comment:} #1}}

\renewcommand{\mod}{\;\operatorname{mod}}
\newcommand{\ord}{\operatorname{ord}}
\newcommand{\TT}{\mathbb{T}}
\renewcommand{\i}{{\mathrm{i}}}
\renewcommand{\d}{{\mathrm{d}}}
\renewcommand{\^}{\widehat}
\newcommand{\HH}{\mathbb H}
\newcommand{\Vol}{\operatorname{vol}}
\newcommand{\area}{\operatorname{area}}
\newcommand{\tr}{\operatorname{tr}}
\newcommand{\mult}{\mu} %weighted multiplicity
\newcommand{\amp}{\beta_a} % \mult/2\sinh(\log NP}
\newcommand{\norm}{\mathcal N} % norm =(\frac{ n+\sqrt{n^2-4}} 2)^2
\newcommand{\LH}{L} % inverse width of window
\newcommand{\Wos}{S}%{S_{a,f,\LH}}
\newcommand{\OPN}{\operatorname{Op}_N}
\newcommand{\HN}{\mathcal H_N}
\newcommand{\UN}{U_N}
\newcommand{\CN}{C_N} %Hecke ops
\newcommand{\intinf}{\int_{-\infty}^\infty}
\newcommand{\ave}[1]{\left\langle#1\right\rangle} %  average
\newcommand{\Var}{\operatorname{Var}}
\newcommand{\sym}{\operatorname{Sym}}
\newcommand{\disc}{\operatorname{disc}}
\newcommand{\CA}{{\mathcal C}_A}
\newcommand{\cond}{\operatorname{cond}} % conductor

\title[Fluctuations in short windows]
{Matrix elements for the quantum cat map: 
Fluctuations in short windows}
\author{P\"ar Kurlberg, Lior Rosenzweig  and Ze\'ev Rudnick}

\address{Department of Mathematics, Royal Institute of Technology,
SE-100 44 Stockholm, Sweden}
\email{kurlberg@math.kth.se}

\address{Raymond and Beverly Sackler School of Mathematical Sciences,
Tel Aviv University, Tel Aviv 69978, Israel}
\email{rosenzwe@post.tau.ac.il}

\address{Raymond and Beverly Sackler School of Mathematical Sciences,
Tel Aviv University, Tel Aviv 69978, Israel}
\email{rudnick@post.tau.ac.il}

\date{July 2, 2007}
\thanks{P.K. was partially supported by a grant from the G\"oran
Gustafsson Foundation, the Royal Swedish Academy of Sciences, and the
Swedish Research Council. L.R. and Z.R. were supported by THE ISRAEL
SCIENCE FOUNDATION (grant No. 925/06) }  

\begin{abstract}
We study fluctuations of the matrix coefficients for the quantized
cat map. We consider the sum of matrix coefficients corresponding to
eigenstates whose eigenphases lie in a randomly chosen window,   
assuming that the  length of the window shrinks with Planck's constant. 
We show that if the length of the window is smaller than
the square root of Planck's constant, but larger than the separation
between distinct eigenphases, then 
the variance of this sum is proportional to
the length of the window, with a proportionality constant which
coincides with the variance of the individual matrix elements
corresponding to Hecke eigenfunctions. 

\end{abstract}
\maketitle

%{\bf ***Extra References:} \cite{barnett06} \cite{luo-sarnak-quantum-variance}
%\cite{feingold-peres-matrix-elements} \cite{keating-cat-maps-91}

\section{Introduction}
\subsection{Background} 
Much effort has been expanded in recent years to study 
quantum wave functions of  classically chaotic systems in the
semiclassical limit. 
One well known result is that the matrix elements of 
smooth observables concentrate around the classical average of the observable, 
at least in the mean square \cite{Sch, CdV, Zelditch}; this is known as
the ``Quantum Ergodicity theorem'' and is valid in great generality. 
A harder problem, known as ``Quantum Unique Ergodicity'' (QUE), is the
question whether {\em all} matrix elements converge 
to the classical average of the observable. This is expected to hold
for any negatively curved surfaces \cite{RS}, but unlike the case of
Quantum Ergodicity, there are no general results available here. 
The  only rigorous results available concern special arithmetic systems,  
namely cat maps and some special compact surfaces
of constant negative curvature, uniformized by unit groups of rational
quaternion algebras. In these cases  many quantum symmetries exist, and QUE  is
now known for eigenfunctions of the desymmetrized system, 
\cite{KR1,Lindenstrauss}. The complexity of  the problem  increases as
we increase the number of degrees of freedom and Kelmer \cite{Kelmer}
found systematic deviations from QUE for higher dimensional cat maps. 
Without incorporating the symmetries, QUE is violated for the
two-dimensional cat map \cite{FNB}.  

An important problem is to understand the rate of convergence to the
classical average. It has been suggested by Feingold and Peres
\cite{FP} that for generic systems with $D$ degrees of freedom, the
variance of the matrix elements about their mean decays with Planck's
constant $\hbar$ as  $\hbar^D$, 
with a prefactor given in terms of the  autocorrelation function of
the classical observable. Several (non-rigorous) arguments where given
for this by Eckhardt et al \cite{EFKAMM}. 
%which strongly inspired our present work. 
For an extensive numerical test of the Feingold-Peres conjecture, 
see Barnett \cite{barnett06}. %, to which we refer for further discussion. 
Rigorous results towards this conjecture are only available for
arithmetic systems - the modular domain \cite{LS, Zhao} and the cat map
\cite{KR-annals, Kelmer}. In both cases arithmetic deviations
from the conjecture are found.  

Once one knows the variance, it is natural to believe that the 
{\em normalized} matrix elements fluctuate randomly about the mean. 
In this paper, we study the fluctuations of matrix elements of the
quantized cat map, 
%a model quantum system with chaotic classical analogue, 
by studying the variance of  {\em sums} of the matrix coefficients
over randomly chosen energy windows. Our findings is that indeed there
is considerable cancellations in these sums,  
%of matrix elements over suitable energy windows, 
consistent with a supposition that the 
signs of the normalized matrix elements behave randomly. 
We will describe the results in detail once we recall the model.

% We study ``partial traces'', or ``traces in short intervals'' of
% matrix elements for the quantized cat map.

%In this paper, we study the fluctuations of matrix elements of a model
%quantum system with chaotic classical analogue, by studying the
%variance of  sums of
%the matrix coefficients over randomly chosen energy windows.

%The model system is the quantized cat map:

\subsection{The quantum cat map} 
The quantized cat map is a model quantum system with chaotic classical
analogue, first investigated by Hannay and Berry \cite{HB} and studied
extensively since, see e.g. \cite{Keating 91, DEGI, KR1, FNB,
Rudnick-montreal}.  
While the classical system displays generic chaotic properties, the
quantum system behaves non-generically in several aspects, such as the
statistics of the eigenphases, and the value distribution of the
eigenfunctions \cite{KR-imrn}.  

We review some of the details of the system in a form suitable for
our purposes, see e.g. \cite{DEGI, KR1, Rudnick-montreal}.  
Let $A$ be a linear hyperbolic toral automorphism, that is,
$A\in SL_2(\Z)$ is an integer unimodular matrix with distinct real
eigenvalues. We assume $A\equiv I\mod 2$. 
Iterating the action of $A$ on the torus $\TT=\R^2/\Z^2$ gives  a
dynamical system, which is highly chaotic.  
The quantum mechanical system includes an integer $N\geq 1$, the inverse
Planck constant, (which we will take to be prime), 
an $N$-dimensional state space $\HN\simeq L^2(\Z/N\Z)$, and a unitary map  
$U=\UN(A)$ of $\HN$, which is the quantization $A$. 
Fix a smooth real-valued  observable $f\in C^\infty(\TT)$, 
%which in fact we will take to be  a trigonometric polynomial (that is
%having only finitely many Fourier modes), 
which we will assume has zero mean: 
$
\int_{\TT} f(x)dx =0\;,
$ 
and let $\OPN(f)$ be its
quantization, which is a self-adjoint operator on $\HN$.
%(which we will in fact assume is a trigonometric polynomial?)
Let $\psi_j$ be an orthonormal basis of
eigenstates for $U$ with eigenvalues $e^{2\pi i\theta_j}$:
$U\psi_j=e^{2\pi i\theta_j}\psi_j$, and $\langle \OPN(f) \psi_j,\psi_j
\rangle$ the (diagonal) matrix elements.
%The matrix elements are real numbers if the observable $f$ is real-valued.

Let  $\ord(A,N)$ be the
least integer $r\geq 1$ for which $A^r\equiv I \mod N$.
When $N$ is prime the distinct eigenphases $\theta_j$ 
are evenly spaced (with at most one exception)
with spacing $1/\ord(A,N)$, and in fact, 
%after removing  a phase in the quantization of $U_N(A)$, 
the distinct eigenphases are all of the form $j/\ord(A,N)$. 
The eigenspaces all have the same dimension (again with at most one
 exception) which is $(N\pm 1)/\ord(A,N)$.

For fixed small $\epsilon>0$, as $N\to\infty$ through a sequence of
values such that  $\ord(A,N)>N^\epsilon$ (which is valid for most
values of $N$; in fact $\ord(A,N)>N^{1/2+o(1)}$ for almost all $N$, 
c.f. \cite[Lemma 15]{KR2}),
all the matrix elements converge to the phase space average
$\int_{\TT} f(x)dx$ of the observable $f$ \cite{KR2, Bourgain-cat}. 
%It is important to keep in mind that this is not always be the case: 
(However, note that there are ``scars'' found for values of $N$ where
$\ord(A,N)$ is 
logarithmic in $N$, see \cite{FNB}.)

To study the fluctuations of the matrix elements, 
%about $\int_{\TT} f(x)dx=0$,
we study the sum of diagonal matrix elements of
$\OPN(f)$ over eigenphases lying in a random window of length
$1/L$ around $\theta$. 
%$$P(\theta) = \sum_{|\theta_j-\theta|<\frac 1{2L} }
%\langle \OPN(f)\psi_j,\psi_j \rangle  \;.$$
More generally we consider a window function,
constructed by taking a fixed non-negative and even function $h \in
L^2([-\frac 12,\frac 12])$  
and setting $h_L(\theta) :=\sum_{m\in \Z}
h(L(\theta-m))$, which is periodic and localized in an
interval of length $1/L$. We further normalize so that $\intinf h(x)^2dx=1$,  
and hence $\int_0^1 h_L(\theta)^2 d\theta %= \frac 1L \intinf h(x)^2 dx
= 1/L$. 
Then set 
\begin{equation}\label{defn P}
P(\theta) := \sum_{j=1}^N h_L(\theta-\theta_j)
\langle \OPN(f)\psi_j,\psi_j \rangle  \;.
\end{equation}
Note that $P(\theta)$ is independent of choice of basis.
If $f=g-g\circ A$ (where $g\circ A$ is defined by $(g\circ
A)(x)=g(Ax)$) is a cocycle then, as follows from ``Egorov's 
theorem'', all matrix elements 
$\langle\OPN(f)\psi_j,\psi_j \rangle=0 $ vanish  and so
$P(\theta)\equiv 0$ in this case.  

%If in addition the  matrix elements randomly about the mean of $f$,
%then we  expect  the typical size of
%$P(\theta)-\frac NL\int_{\TT}f(x)dx$ to be $1/\sqrt{L}$.
%\begin{verbatim}
%this is not true if f has constant term, 
%since then fluctuations will be O(1) ??
%\end{verbatim}
%To test this, 

The expected value of $P(\theta)$, when we pick $\theta$
randomly and uniformly in the unit interval, is $\int_0^1
P(\theta)d\theta = 0$, since we  assume that $\int_{\TT}
f(x)dx=0$.  (Strictly speaking, for $f$ smooth we only have $\int_0^1
P(\theta)d\theta = O_{f,R}(N^{-R})$ for all $R>0$, see
section~\ref{subsec:formula-variance} for more details.)

%from the definition  as 
%\begin{equation*}
%\begin{split}
%\int_0^1 P(\theta)d\theta &=\frac{\intinf h(x)dx}{L} \sum_{j=1}^N 
%\langle \OPN(f)\psi_j,\psi_j \rangle  \\
%&= 
%\frac{1}{L}\intinf h(x)dx \left(N\int_{\TT} f(x)dx +O(N^{-\infty})
%\right)
%\end{split}
%\end{equation*}

We will study the variance  of $P(\theta)$.  
To describe it, %the asymptotic variance $\Var(P)$,
we introduce the quadratic form associated to the
matrix $A$ by 
%$A=\begin{pmatrix}   a&b\\c&d\end{pmatrix}$    
\begin{equation*}
%\label{quad form}
Q(x)=\omega(x,xA) %cx^2+(d-a)xy-by^2 \;.
\end{equation*}
where $\omega(x,y) = x_1y_2-x_2y_1$ is the standard symplectic form. 
If the Fourier expansion\footnote{With $k=(k_1,k_2)$, $\^f(k) :=
  \int_{0}^1 \int_{0}^1 f((x_1,x_2)) e^{-2 \pi i(k_1x_1 + k_2x_2)}dx_1dx_2$.}
of the observable is 
$f(x) = \sum_{k\in \Z^2} \^f(k) e(kx)$,
where in what follows we abbreviate $e(z):=e^{2\pi i z}$, 
set
% $$f^{\#}(\kappa):=\sum_{k:Q(k)=\kappa} (-1)^{k_1 k_2} \^f(k)$$
%and
\begin{equation}\label{def Carith} 
C_{arith}(f):=\sum_{k,k'\in \Z^2}\sum_{Q(k)=Q(k')} (-1)^{k_1 k_2+k'_1 k'_2} \^f(k)\overline{\^f(k')} \;.
%\sum_{\kappa\neq 0} f^{\#}(\kappa)^2 \;.
\end{equation}

\subsection{Results} 
We can now formulate our main result:
\begin{thm}\label{main thm}
Fix $f\in C^\infty(\TT^2)$ of zero mean. 
Consider any sequence of primes $N$ for which 
$\ord(A,N)/\sqrt{N}\to\infty$. Assume that  $L<2\ord(A,N)$.  
Then as  $N\to \infty$,
$$
\Var(P)= \frac 1L C_{\mbox{arith}}(f) +O \left( \frac{\sqrt{N}}{L^2} \right)\;.
$$
\end{thm}

It is easy to check that $C_{arith}(f)$ vanishes if $f=g-g\circ A$ is
a cocycle.  In that case $P=0$ so Theorem~\ref{main thm} has no
content. 
In case $C_{arith}(f)\neq 0$, if we further assume
$L/\sqrt{N}\to \infty$   we get 
\begin{cor}
Under the assumptions of Theorem~\ref{main thm}, if 
$C_{arith}(f)\neq 0$, and $L/\sqrt{N}\to \infty$ then
$$\Var(P) \sim \frac{C_{arith}(f)}{L} \;.$$ 
\end{cor}

The prefactor 
$C_{arith}(f)$ coincides with the asymptotic variance of the
normalized matrix elements of $\OPN(f)$ when computed in the Hecke
basis \cite{KR-annals}. 
To explain this, note that if $1/L$ is smaller then
the minimal separation between distinct eigenphases, that is if
$L>\ord(A,N)$, then from the definition \eqref{defn P} we get 
\begin{equation}\label{spectral P}
\int_0^1 P(\theta)^2 d\theta = \frac 1L 
\sum_{j=1}^{\ord(A,N)} \left|\sum_{k : \theta_k=\frac {j}{\ord(A,N)} } 
\langle \OPN(f)\psi_k,\psi_k \rangle  \right|^2 \;,
\end{equation}
the inner sum being over an orthonormal basis of the eigenspace
corresponding to a given eigenphase. 
In particular, once $L$ exceeds $\ord(A,N)$, the dependence of
$P(\theta)$ on $L$ is essentially trivial, and thus we may (and shall)
restrict $L<2\ord(A,N)$. 
If in addition we have $\ord(A,N)=N\pm 1$, then (almost) all the inner
sums in \eqref{spectral P} collapse to a single term, and thus we find that
when $\ord(A,N)=N\pm 1$ is maximal and $L>N+1$ that 
$$
\int_0^1 P(\theta)^2 d\theta  = \frac 1L \sum_{j=1}^N 
\left| \langle \OPN(f)\psi_j,\psi_j \rangle  \right|^2 +O(\frac 1{NL})\;.
$$
Thus we recover the variance of the individual matrix elements (a
trick  used by Berry \cite{Berry}), which in turn was shown to be
$C_{arith}(f)$ in \cite{KR-annals}.

For comparison of our results with those expected of {\em generic} quantum
chaotic systems,  consider the case where we take $A$ to be a
non-linear perturbation of a cat map, %by a nonlinear shear, 
the perturbation sufficiently small so that the map remains hyperbolic 
\cite{de Matos and de Almeida}. The resulting system
is expected to have generic spectral statistics (depending on the
symmetries of the map) 
\cite{de Matos and de Almeida, Keating-Mezzadri}. 
Arguing as in \cite{EFKAMM, dCKR} one then expects that 
the variance of $P(\theta)$ in this case 
is asymptotic to $\frac 1L C_{gen}(f)$, where
\begin{equation}\label{Cgen}
C_{gen}(f)=\sum_{t=-\infty}^\infty \int_{\TT} (f\circ A^t)(x)f(x)dx
\end{equation}
is the classical autocorrelation function.  
%($ f_0=f-\int_{\TT}f(x)dx$).
%This is also expected to be the asymptotic variance for the diagonal
%matrix elements of $\OPN(f)$ in the generic case.
To compare with the arithmetic variance $C_{arith}(f)$ 
in \eqref{def Carith}, note that when 
$A$ is {\em linear}, we may write \eqref{Cgen} 
in terms of the Fourier expansion of the observable $f$ as 
$$
C_{gen}(f)=\sum_{k,k'\in \Z^2}\sum_{k'\sim k} (-1)^{k_1 k_2+k'_1
k'_2} \^f(k)\overline{\^f(k')} \;,
$$
where the sum $k'\sim k$ is over pairs of frequencies which lie in
the same $A$-orbit: $k'=kA^t$ for some integer $t$. This condition 
implies the condition $Q(k')=Q(k)$ which enters in the sum \eqref{def Carith} 
for the arithmetic variance $C_{arith}(f)$ (and also implies that 
$(-1)^{k_1 k_2+k'_1 k'_2}=1$).  
%While we have no tools to study the generic case, it would be of
%interest to study the analogous problem for the modular domain.

\subsection{About the proof} 
As we explain in section~\ref{sec: proof of thm}, the variance of $P$
can be written as 
\begin{equation}\label{varexp0}
\Var(P) =
%\ave{|P(\theta)-\frac 1L \^h(0) \tr \{ \OPN(f)|^2}\} =
\frac 1{L^2}\sum_{t\in \Z-\{0\}} \^h\left(\frac tL\right)^2
|\tr\{\OPN(f)U^{-t}\}|^2 \;.
\end{equation}

One approach to evaluating \eqref{varexp0}, used in \cite{EFKAMM}, is
to use a trace formula expressing $\tr\{\OPN(f)U^t\}$ as a sum over
periodic orbits of the map $A$ with certain phases, where the number
of summands grows exponentially in $t$.  
%(Also see \cite{Keating 91, Keppeler-Marklof-Mezzadri} for a trace formula
%specifically for quantized maps of the torus.)
This gives $\Var(P)$ as a sum over pairs of periodic orbits. The
averaging over $\theta$ produces a sum over ``diagonal'' pairs where
the phases cancel,  and the remaining pairs. The obvious diagonal
pairs consist of equal orbits and give the generic answer
$C_{gen}(f)$. To reproduce the correct answer $C_{arith}(f)$ in this
case requires identifying another diagonal family and showing that the
contribution of the remaining pairs is negligible. We have not been
able to do that. 

Instead, we use a different formula for $\tr \OPN(f)U^{-t}$, based on a
formula for the quantum propagator $U_N(A)$ introduced by Kelmer
\cite{Kelmer}   and a expansion in Fourier modes of $f$ to
rewrite \eqref{varexp0} as a double sum over Fourier modes
\begin{equation*}%\label{first sum for P} 
\Var(P) \sim \frac 1{L^2}\sum_{k,k'}
(-1)^{k_1k_2+k'_1k'_2}\^f(k)\overline{\^f(k')} S(k,k') %+O(N^{-\infty}) 
\end{equation*}
where $S(k,k')$ is a certain incomplete exponential sum, which is
trivial for pairs of frequencies with $Q(k)=Q(k')$. 
These pairs of frequencies give the main term of
$C_{arith}(f)/L$ ; this is our new
diagonal approximation, ``dual''  in a sense to the standard one using
periodic orbits. 
To handle the off-diagonal terms, it suffices to
give a non-trivial bound for the exponential sum $S(k,k')$ when
$Q(k)\neq Q(k')$. Using a
standard completion technique, we reduce it to giving a bound for a
certain complete exponential sum. When $N$ is a ``split'' prime, that
is if $A$ is diagonalizable modulo $N$, the required bound is a 
standard result of the Riemann Hypothesis for function fields (proved
by Weil). For the remaining ``inert'' primes, the required bound was
recently established by Gurevich and Hadani \cite{GH}. 
%who invoke the full force of Deligne's Weil II paper \cite{Deligne Weil 2} to
%give the required bound. 
In the appendix, we will give a different proof that only requires
Weil's original methods \cite{Weil48}.

{\bf Acknowledgment:}  We thank Dubi Kelmer for several helpful
comments, and the referees for several suggestions for improving the
presentation. 
 
\section{Prerequisites on cat maps and their quantization}
\subsection{The quadratic form associated to $A$}
Let $A=\begin{pmatrix}a&b\\c&d\end{pmatrix}\in SL_2(\Z)$ be
hyperbolic, and assume that $A\equiv I\mod 2$. Then $A$ preserves the
standard 
symplectic form
$$
\omega(x,y):=x_1y_2-x_2y_1  \;
$$
and thus  $A$ also preserves the quadratic form
\begin{equation}\label{def Q}
Q(x):=\omega(x,xA) = bx_1^2+(d-a)x_1x_2 -c x_2^2 \;,
\end{equation}
which has discriminant $\disc(Q)=(\tr A)^2-4$ (which is even since
$A\equiv I\mod 2$).

\begin{lem}\label{lem: one-dim space}
Let $N$ be an odd prime, and let $A\in SL_2(\Z)$ so that 
$(\tr A)^2-4\neq 0\mod N$. 
Then the space of binary quadratic forms preserved by $A$ is one dimensional.
\end{lem}
\begin{proof}
Passing if necessary to a quadratic extension $\FF$ of the base
field $\Z/N\Z$ over which $A$ is diagonalizable, as we may by our
assumption on $N$, consider the action of $A$ on $2\times2$ matrices 
over $\FF$ via 
$M\mapsto A^TMA$ (where $A^T$ is the transpose of $A$). 
The decomposition of the space of matrices to a direct sum
of the one-dimensional space of skew-symmetric and the
three-dimensional space $\sym^2$ of symmetric matrices (identified
with quadratic forms) is preserved by $A$.

Let $\lambda^\pm $ be the eigenvalues of $A$. Since $N$ is coprime to
$\tr(A)^2-4$,  $A\neq \pm I\mod N$ and in particular $\lambda\neq \pm 1$. 
Let $v_\pm\neq 0$ be the corresponding  eigenvectors: $v_\pm
A=\lambda^\pm v_\pm$ (note that our vectors are row vectors). 
Then the matrices $v_\pm^T v_\pm$ are  four
eigenvectors for the action of $A$ on $2\times 2$ matrices, with
eigenvalues $\lambda^2, 1, 1, \lambda^{-2}$. The skew symmetric
matrices are fixed by $A$ and hence the eigenvalues of $A$ on
$\sym^2$ are $\lambda^2,1,\lambda^{-2}$. The $1$-eigenspace
corresponds to binary quadratic forms which are preserved by $A$.
Since $\lambda\neq \pm 1$,
%we are assuming that $A\neq \pm I\mod N$, we have $\lambda^2\neq 1$
%and so 
we find that the $1$-eigenspace is one-dimensional, proving the claim.
\end{proof}
Thus we find that if $N$ is  coprime to $\disc(Q)$, then any binary
quadratic form preserved by $A$ is a multiple of $Q$ modulo $N$.

\subsection{} 
Let $N$ be a prime not dividing $\disc(Q)=(\tr A)^2-4$. 
% so that $A\neq \pm I\mod N$ 
In \cite{KR1} we described a commutative algebraic group
$\CA(N)\subset SL_2(\Z/N\Z)$, containing $A$, which in the case at
hand is the centralizer  of $A$ in $SL_2(\Z/N\Z)$ 
and coincides with 
the special orthogonal group of the quadratic form
$Q$ given in \eqref{def Q} over the field $\Z/N\Z$. 
The group $\CA(N)$ is isomorphic to
either the multiplicative group of the field $(\Z/N\Z)^*$ (the
``split'' case) or the norm-one elements in a quadratic extension of
$\Z/N\Z$ (the ``inert'' case). Thus $\CA(N)$ has order $N-1$ or $N+1$,
respectively.  Note that if $g\in \CA(N)$ and $g\neq 1$, then the
matrix $g-1$ is 
invertible. 

\subsection{}\label{sec q_t}
As examples of quadratic forms over $\Z/N\Z$ preserved by $A$, consider
for $g\in \CA(N)$, $g\neq 1$, 
$$
q(x;g):= \omega(x(g-1)^{-1},x(g-1)^{-1}g).
$$
Note that $q(\bullet,g)=0$ if $g=-I\mod N$. 
\begin{lem}\label{lem: q multiple of Q}
If $g\neq \pm I \mod N$ then $q(\bullet;g)$ is a nonzero multiple of $Q$.
\end{lem}
\begin{proof}
By Lemma~\ref{lem: one-dim space} it suffices to show that $q(\bullet;g)$ is
preserved by $A$ and is nonzero.
It is preserved
by $A$ since $A$ preserves $\omega$ and commutes with both $g$ and
$(g-I)^{-1}$.  By Lemma~\ref{lem: one-dim space}, $q(\bullet;g)$ is thus a
multiple of the form $Q$. We claim that the multiple is
nonzero mod $N$. To see this, it suffices see that $q(\bullet;g)$ is not
identically zero. If this were the case, then since $\omega$ is
non-degenerate and $g-I$ invertible, we would have that
$yg$ is a scalar multiple of $y$ for all vectors $y$, which
necessarily forces $g$ to be a scalar matrix. Since $\det g=1$, we
must therefore have $g=\pm I\mod N$, contradicting our assumption.
\end{proof}

\subsection{Computing $Q(x)$ and $q(x;g)$} 
 Choose a generator $g_0$ of the
cyclic group $\CA(N)$. 
Passing if necessary to a quadratic extension
of $\Z/N\Z$, write $x=x_+ + x_-$ where $x_\pm $ are
eigenvectors of $g_0$ hence of  $g$ and of $A$, with 
$x_\pm g = \lambda^{\pm 1}x_\pm$, $x_\pm A= \lambda_A^{\pm 1}x_\pm$
($\lambda_A\neq \pm 1$ if $N$ is coprime to $\disc(Q)=(\tr a)^2-4$). 
Then
\begin{equation*}
\begin{split}
Q(x) &= \omega(x,xA) = \omega(x_+ + x_-, \lambda_A x_+ + \lambda_A^{-1} x_-)\\
&= (\lambda_A-\lambda_A^{-1}) \omega(x_-,x_+)\;.
\end{split}
\end{equation*}
Likewise,
\begin{equation*}
\begin{split}
q(x;g) &= \omega (x(g-I)^{-1},x(g-I)^{-1} g) \\
&=
\omega( \frac 1{\lambda-1} x_+ + \frac 1{\lambda^{-1}-1} x_-,
\frac {\lambda}{\lambda-1} x_+
+ \frac{\lambda^{-1}}{\lambda^{-1}-1} x_-  )\\
&=  \frac{\lambda-\lambda^{-1}}
{(\lambda-1)(\lambda^{-1}-1)} \omega(x_-,x_+) \;.
\end{split}
\end{equation*}
In particular we find
\begin{equation}\label{formula for q}
q(x;g) =Q(x) \frac{\lambda-\lambda^{-1}}{\lambda_A-\lambda_A^{-1}}
\frac 1{(\lambda-1)(\lambda^{-1}-1)}=
\frac{Q(x)}{\lambda_A-\lambda_A^{-1}}
\frac{1+\lambda}{1-\lambda}.
\end{equation}

\subsection{Quantum mechanics on the torus}
We recall the basic facts of quantum mechanics on the torus which we
need in the paper,  see
\cite{KR1, Rudnick-montreal} for further details.
Planck's constant is restricted to be
an inverse integer $1/N$, and the Hilbert space of states $\HN$ is
$N$-dimensional, which is identified with $L^2(\Z/N\Z)$ 
with the inner product given by
\begin{equation*}
\langle \phi,\psi \rangle
 := \frac1N \sum_{Q\bmod N} \phi(Q) \, \overline\psi(Q) \;.
\end{equation*}

Classical observables, that is real-valued functions $f\in
C^\infty(\TT)$, give rise
to quantum observables, that is  self-adjoint operators $\OPN(f)$
on $\HN$.
To define these, one starts with translation operators: for 
$n=(n_1,n_2)\in\Z^2$ let $T_N(n)$ be the unitary operator on $\HN$ 
whose action on a wave-function $\psi\in \HN$ is 
\begin{equation*} %\label{action of T(n)}
T_N(n)\psi(Q) = e^{\frac {i\pi n_1 n_2}N}
e\left(\frac{n_2Q}N\right)\psi(Q+n_1) \;.
\end{equation*}
For any smooth function $f\in C^\infty(\TT)$,  define $\OPN(f)$ by 
$$
\OPN(f) := \sum_{n\in\Z^2} \widehat f(n) T_N(n)
$$
where $\widehat f(n)$ are the Fourier coefficients of $f$.
%If $f$ is a trigonometric polynomial, then for $N>N_0(f)$
%sufficiently large 
The trace of $\OPN(f)$ is 
\begin{equation}\label{eq: trace OPN(f)} 
\tr \{ \OPN(f) \} = N\int_{\TT} f(x)dx +O_f ( N^{-\infty}) 
\end{equation}
where the term $O_f(N^{-\infty})$ is one that is bounded by $N^{-R}$
for any $R>0$, the implied constant depending on $f$ and $R$. 

\subsection{A  formula for the  quantum propagator}
\label{sec:kelmers-formula}
For any $B\in SL_2(\Z)$, $B\equiv I\mod 2$, the quantum propagator
$U_N(B)$ is a unitary map of $\HN$ satisfying Egorov's formula
$$U_N(B)^*\OPN(f)U_N(B)=\OPN(f\circ B)$$
for all observables $f\in C^\infty(\TT)$.
This property defines the propagator only
up to a phase, which will be of no interest to us.
%since it plays no role in the formula \eqref{}  for the variance.

A useful formula for the propagator, known in the context of the Weil
representation (cf. \cite{Moeglin}), and introduced for cat maps by
Kelmer \cite{Kelmer}, is the following:
for any $B\in SL_2(\Z)$, and $N$ odd, the quantum propagator
$U_N(B)$ is given by
\begin{equation*}
\label{eq:kelmer}
U_N(B)=\frac1{N |\ker_N(B-I)|^{1/2}}
\sum_{n\in(\Z/N\Z)^2}e \left(\frac{\omega(n,n B)}{2N}\right)T_N(n(I-B))
\end{equation*}
where $\ker_N(B-I)$ denotes the kernel
of the map $B-I$ on $\Z^2/N\Z^2$.
We apply this when  $N$ is a prime 
not dividing\footnote{Since $A\equiv I\mod 2$, $\disc(Q)$ is even and 
hence such $N$ is odd.}   $\disc(Q)$,  
$B=A^t$ (where $A^t$ is the $t$-th power of $A$)  so that $A^t\neq  I\mod N$. 
Note that $|\ker_N(A^t-I)|=1$ since in the group $\CA(N)$, 
if $g\neq 1$, then the matrix $g-1$ is invertible. 
%by Lemma~\ref{lem:invertible}.
Thus %, up to an unimportant phase,
\begin{equation}\label{def:cat_map_quantization}
U_N(A^t) = \frac 1N\sum_{n\in (\Z/N\Z)^2}
e \left(\frac{\omega(n,n A^t)}{2N}\right)T_N(n(I-A^t))
\end{equation}

\begin{lem}\label{Kelmer's lemma}
Let  $A\in SL_2(\Z)$ be hyperbolic, and assume that $A\equiv I\mod 2$.
Then for any prime $N$ not dividing $\disc(Q)$ and integer $t$ 
such  that $A^t\neq I\mod N$,  
we have
\begin{equation*}\label{kelmer's formula}
\tr \{ \OPN(f)U_N(A^t) \} = \sum_k (-1)^{k_1k_2} \^f(k)
e \left(\frac{\overline{2}q(k;A^t)}{N} \right)
\end{equation*}
%where $e(\phi(t))$ is a phase factor that comes from a choice of
%quantization of $A^t$,
where  $\overline{2}$ is the inverse of $2\mod N$. 
%and $n=n(t,k)$ is the unique solution modulo $N$ of
%$k\equiv n(A^t-I)\mod N$.
\end{lem}

%\begin{lem}\label{lem:quant_trace}
%\end{lem}
% \begin{remark}\label{rem:even}
% Since $A\equiv I\pmod2$ then $\omega(n,nA^t)\equiv0\pmod2$, and by
% the Chinese Remainder Theorem, we have that $e(\frac
% x{2N})=e(\frac x2)e(\frac{\bar2x}N)$, then we can consider
% $e(\frac{\omega(n,nA^t)}{2N})=e(\frac{\bar2\omega(n,nA^t)}{N})$
% \end{remark}

\begin{proof}
It suffices to show that
\begin{equation}\label{eq:quant_trace}
\tr\{T_N(k)U(A^t)\}=(-1)^{k_1k_2}e\left(\frac{\overline{2}\omega(n,nA^t)}{N}\right)
\end{equation}
where $n$ is such that $k=n(A^t-I)\pmod N$.
Using   \eqref{def:cat_map_quantization}
% $$ U_N(A^t)=\frac1{N} \sum_{n\in(\Z/N\Z)^2}
% e(\frac{\omega(n,nA^t)}{2N})T_N(n(I-A^t)) $$
we find
$$
\tr \{T_N(k)U_N(A^t)\} =\frac1{N} \sum_{n\in(\Z/N\Z)^2}
e\left(\frac{\omega(n,nA^t)}{2N}\right) \tr \{ T_N(k)T_N(n(I-A^t)) \} \;.
$$
As is easy to see from the definition (see  Lemma~4 and (2.6) in \cite{KR1}),
$$
\tr \{T_N(n)T_N(m) \}=
\begin{cases}
     (-1)^{m_1m_2+n_1n_2}N & \text{if $n\equiv -m\pmod N$,}\\
       0 & {\rm{otherwise.}}
   \end{cases} 
$$
%using this with $n=k,m=-n(I-A^t)$, we get that
Thus $\tr \{T_N(k)T_N(n(I-A^t))\}=(-1)^{k_1k_2}N$ if $k=-n(I-A^t)\pmod N$
%(since $A\equiv I\pmod2$, we have that $-n(I-A^t)\equiv0\pmod2$),
and 0 otherwise.
Now  if $A^t\neq I\mod N$, such $n$ as above
exists and is unique since $A^t-I$ is invertible.  
%by Lemma~\ref{lem:invertible}.  
Therefore
\begin{equation}\label{eq:quant_trace2}
\tr\{T_N(k)U(A^t)\}=(-1)^{k_1k_2}e\left(\frac{\omega(n,nA^t)}{2N}\right) \;.
\end{equation}
Now if $b$ is an even integer and $N$ is odd then
$e(\frac{b}{2N})=e(\frac{\overline{2}b}{N})$.
Applying this to \eqref{eq:quant_trace2}  with
$b=\omega(n,nA^t)$ which is even since
$A\equiv I\mod 2$, we end up with formula \eqref{eq:quant_trace}.
\end{proof}

%\begin{remark}\label{rem:Hecke orbit}
%According to Lemma~\ref{lem: q multiple of Q}, if $A^t\neq \pm I\mod
%N$ then 
%$$e(\frac{\omega(n,nA^t)}{2N})=e(\frac{\omega(n',n'A^t)}{2N})$$
%for all $k'$ satisfying $k'=Bk$ for some $B\in SL_2(\Z/N\Z)$ that
%commutes with $A$ in $SL_2(\Z/N\Z)$, and this condition is
%satisfied
%if and only if %$\omega(k,kA)=\omega(k',k'A)\mod N$, that is iff
%$Q(k)\equiv Q(k')\mod N$. By \cite[?????]{KR-annals}, this is equivalent
%to $k$ and $k'$ lying in the same Hecke orbit.
%\end{remark}

\section{Proof of Theorem~\ref{main thm}}\label{sec: proof of thm}

\subsection{A formula for the variance} 
\label{subsec:formula-variance}
%It is convenient to use a smooth count:
Fix a non-negative, even, 
%smooth\footnote{This is not necessary; all our
%conclusions remain valid when we take $h$ %=\mathbf 1_[-1/2,1/2]$ 
%to be an indicator function.} 
test function $h$, 
supported in $[-\frac 12,\frac 12]$  
and normalized so that  $\intinf h(x)^2dx=1$. 
Set
$$
h_L(x):=\sum_{k\in \Z} h(L(x-k))
$$
which is then a periodic function, localized on the scale of
$1/L$, and  $\int_0^1 h_L(\theta)^2d\theta =1/L$.    
The Fourier expansion of $h_L$ is (in $L^2$ sense)
$$
h_L(x)= \frac 1L\sum_{t\in \Z} \^h \left(\frac tL \right) e(tx)\;.
$$
where $\^h(y) = \int_{-\infty}^{\infty} h(x) e(-xy) \, dx$.  

Let $N$ be a prime which does not divide $\disc(Q)= (\tr A)^2-4$.  
Let
$$
P(\theta) := \sum_j h_L(\theta-\theta_j) \langle \OPN(f)\psi_j,\psi_j \rangle
$$
which is a  sum of matrix elements on a window of size
$1/L$ around $\theta$.
Then, in $L^2$ sense, and with $U=U_N(A)$, we have
\begin{equation}
P(\theta)  = \frac 1L\sum_{t\in \Z} e(t\theta)\^h\left(\frac tL\right)
\tr \{\OPN(f)U^{-t}\} \;.
\end{equation}
(Note that
$|\tr \{\OPN(f)U^{-t}\}|$ is uniformly bounded in $t$.)
The mean value of $P(\theta)$  is
%is
$$%\ave{P(\theta)} =
\int_0^1 P(\theta)d\theta = 
\frac 1L \^h(0) \tr \{ \OPN(f) \}
= %\frac {\^h(0)N}{L} \int_{\TT^2} f(x)dx +
O_{f}(N^{-\infty})
%=0
$$
according to \eqref{eq: trace OPN(f)}. 
%since if $f$ is a trigonometric polynomial then for $N>N_0(F)$, 
%$\tr\OPN(f)=N\int_{\TT} f(x)dx = 0$. 
Thus the variance can be written as
\begin{equation}\label{varexp}
\Var(P) =
%\ave{|P(\theta)-\frac 1L \^h(0) \tr \{ \OPN(f)|^2}\} =
%\frac 1{L^2}\sum_{t\in \Z} \^h\left(\frac tL\right)^2
\frac 1{L^2}\sum_{t\in \Z-\{0\}} \^h\left(\frac tL\right)^2
%|\tr\{\OPN(f)U^{-t}\}|^2 + O_{f}(N^{-\infty})\;.
|\tr\{\OPN(f)U^{-t}\}|^2\;.
\end{equation}
%This formula remains valid for a more general class of functions $h$. 
%in particular for the indicator function of an interval.

\subsection{Computing the variance} 
%One approach to evaluating \eqref{varexp}, used in \cite{EFKAMM}, is
%to use a trace formula expressing $\tr\{\OPN(f)U^t\}$ as a sum over
%periodic orbits of the map $A$ with certain phases, where the number
%of summands grows exponentially in $t$.  
%%(Also see \cite{Keating 91, Keppeler-Marklof-Mezzadri} for a trace formula
%%specifically for quantized maps of the torus.)
%This gives $\Var(P)$ as a sum over pairs of periodic orbits. The
%averaging over $\theta$ produces a sum over ``diagonal'' pairs where
%the phases cancel,  and the remaining pairs. The obvious diagonal
%pairs consist of equal orbits and give the generic answer
%$C_{gen}(f)$. To reproduce the correct answer $C_{arith}(f)$ in this
%case requires identifying another diagonal family and showing that the
%contribution of the remaining pairs is negligible. We have not been
%able to do that. 

%Instead, we use the results of Section~\ref{sec:kelmers-formula} as
%follows: 
Let $\ord(A,N)$ be the least integer $r\geq 1$ so that
$A^r\equiv I\mod N$.  It is a divisor of $|\CA(N)|$, that is of either
$N-1$ or $N+1$.  We will rewrite \eqref{varexp} as
\begin{equation}\label{modular sum}
\Var(P) =\frac 1{L^2} \sum_{\tau \mod \ord(A,N)} \Gamma(\tau)  
%|\tr\{\OPN(f)U^{-\tau}\}|^2 + O_{f}(N^{-\infty})\;,
|\tr\{\OPN(f)U^{-\tau}\}|^2\;,
\end{equation}
where 
\begin{equation}\label{defn Gamma}
%\Gamma(\tau) = \sum_{\substack{ t\in \Z\\ t\equiv \tau \mod
\Gamma(\tau) = \sum_{\substack{ t\in \Z -\{0\}\\ t\equiv \tau \mod
\ord(A,N)}}  \^h\left(\frac tL\right)^2 
\;.
\end{equation}
We may omit the term $\tau \equiv 0\mod \ord(A,N)$ from the sum 
\eqref{modular sum} at the cost of introducing an error of
$O(N^{-\infty})$, since then
$\tr \{\OPN(f)U^{-\tau}\} = \tr \{\OPN(f)\}=O(N^{-\infty})$ while
$|\Gamma(\tau)|\ll L$. 
% for $N>N_0(f)$. 

Now  we use Lemma~\ref{Kelmer's lemma} to rewrite $\tr \{
\OPN(f)U_N(A^t)\}$, where we 
replace $U_N(A^t)$  by $U_N(A)^t$ after introducing a phase (which
can be ignored as we are taking absolute values), and replacing
$t$ by $-t$ in \eqref{varexp}, as we may since $h$ is even. 
%: Let $\omega(x,y):=x_1y_2-x_2y_1$ be the standard
% symplectic form, which is preserved by $A$. Then
%Equation \eqref{kelmer's formula}
%. Inserting \eqref{kelmer's formula} into \eqref{varexp} gives
The result is that 
\begin{equation}\label{first sum for P} 
\Var(P) = \frac 1{L^2}\sum_{k,k'} 
(-1)^{k_1k_2+k'_1k'_2}\^f(k)\overline{\^f(k')} 
S(k,k') +O(N^{-\infty}) 
\end{equation}
%with
%\begin{equation}
%S(k,k') = \sum_{\tau\neq 0\mod \ord(A,N)} \Gamma(\tau) 
%e(\frac{\omega(n,nA^\tau)-\omega(n',n'A^\tau)}{2N}) \;.
%\end{equation}
%where for each $k$, $n=n(k,\tau)$ is the unique solution of
%$k=n(A^t-I)\mod N$.
where, in the notation of \S\ref{sec q_t}, we have
\begin{equation}\label{S general case} 
S(k,k') = \sum_{t\neq 0\mod \ord(A,N)}  \Gamma(t) 
e\left(\frac{\overline{2}(q(k;A^t)-q(k';A^t))}{N}\right) \;.
\end{equation}
We have $|S(k,k')|\ll NL$ since $|\Gamma(t)|\ll L$. Thus we may, using
rapid decay of the Fourier coefficients $\^f(k)$, truncate the sum
\eqref{first sum for P}  at frequencies at most $N^{1/4}$ to get 
\begin{equation}\label{truncated sum for P} 
\Var(P) = \frac 1{L^2}\sum_{|k|,|k'|<N^{1/4}}
(-1)^{k_1 k_2+k'_1k'_2}\^f(k)\overline{\^f(k')}S(k,k') +O(N^{-\infty}) \;.
\end{equation}
%This displays  the trace as a finite sum. 

\subsection{Diagonal terms} 
The sum $S(k,k')$ is trivial if the phase difference
%$\omega(n,nA^t)-\omega(n',n'A^t)$ vanishes $\mod 2N$, that is if and
%only if 
$q(k;A^t)-q(k';A^t)$ vanishes $\mod N$ for all $t$. 
By Lemma~\ref{lem: q multiple of Q},
this happens if and only if
%$k$ and $k'$ lie in  the same Hecke orbit, equivalently that
%for the quadratic form $Q$ defined in \eqref{quad form}
we have $Q(k)\equiv Q(k') \mod N$. 
For the frequencies appearing in \eqref{truncated sum for P}, we have
$|Q(k)|,|Q(k')|\ll \sqrt{N}$ by Cauchy-Schwartz, and hence the
congruence $Q(k)\equiv Q(k') \mod N$ forces that 
%Since $f$ is a trigonometric polynomial, 
this latter condition becomes an equality $Q(k)=Q(k')$. 
%for $N>N_0(f)$. 

These ``diagonal'' pairs of frequencies with $Q(k)=Q(k')$ 
give a contribution of
$$
\frac 1{L^2}\sum_{Q(k)=Q(k')}  (-1)^{k_1k_2+k'_1k'_2}\^f(k)\overline{\^f(k')} 
\sum_{t\neq 0\mod \ord(A,N)} \Gamma(t) 
$$
(we may drop the condition $|k|,|k'|<N^{1/4}$ at a cost of
$O(N^{-\infty})$). We claim that 
$$
\sum_{t\neq 0\mod \ord(A,N)} \Gamma(t) = L +O(1) \;.
$$
%as $N,L\to \infty$. ????
To see this, write 
$$
\sum_{\substack{t\mod \ord(A,N)\\ t\neq 0\mod \ord(A,N)}} \Gamma(t)  
= \sum_{t\in \Z} \^h \left(\frac tL \right)^2 - 
\sum_{j\in \Z} \^h \left(\frac {\ord(A,N)}{L}j \right)^2  \;.
$$
Now 
$$
\sum_{t\in \Z} \^h\left(\frac tL \right)^2 = L^2\int_0^1 h_L(\theta)^2d\theta = L 
$$
and 
$$
\sum_{j\in \Z} \^h\left(\frac {\ord(A,N)}{L}j\right)^2 =O(1)
$$
since $L<2\ord(A,N)$. Thus the pairs of frequencies with $Q(k)=Q(k')$
give a total contribution of 
\begin{equation}\label{eq:main term} 
\frac 1L C_{arith}(f) +O\left(\frac 1{L^2} \right) \;.
\end{equation}

\subsection{Off-diagonal terms} 
For the remaining pairs of frequencies, where $Q(k)\neq Q(k')$, 
the sum $S(k,k')$ is a certain incomplete exponential sum.
% whose complete version is akin to a (twisted) Kloosterman sum.

\begin{prop}\label{prop:S is small}
If $Q(k)\neq Q(k')$ then 
$$
|S(k,k')| \ll \sqrt{N} \;.$$
\end{prop}
Assuming we have this,  
the off-diagonal pairs will then contribute at most 
$O(\frac{\sqrt{N}}{L^2})$. 
Thus in combination with \eqref{eq:main term} we get 
$$
\Var(P) = \frac{C_{arith}(f)}{L} + O\left(\frac{\sqrt{N}}{L^2}\right)
$$
which gives Theorem~\ref{main thm}. 

%and if we assume that $C_{arith}(f)\neq 0$ then we are done.
%\begin{verbatim}
%what if C_{arith}(f) = 0 ???
%\end{verbatim}

%The proof of Proposition~\ref{} will be given in the following section.

%\section{Completing exponential sums} 

To prove Proposition~\ref{prop:S is small}, we will need the following
result. 
\begin{lem}\label{lem: bound E(j)} 
Let $k,k'\in \Z^2$, $Q(k)\neq Q(k')$. 
Define %\rmk{It should be $jt$ and not $j\tau$ in 3.8, right!?} 
\begin{equation}\label{defn E(j)}
E_A(j) = \sum_{0\neq t\mod \ord(A,n)} e\left(\frac {jt}{\ord(A,N)}\right)
e\left(\frac{\overline{2} (q(k;A^t)-q(k';A^t))}{N}\right) \;.
\end{equation}
Then for all $j$, 
$$
|E_A(j)| \leq 2\sqrt{N} \;.
$$
\end{lem}
\begin{proof}

For each multiplicative character
$\chi$ of $\CA(N)$, define the complete sum
\begin{equation}\label{def complete sum}
E(\chi):=\sum_{1\neq y\in \CA(N)} \chi(y) 
e\left(\frac{\overline{2}(q(k;y)-q(k';y))}{N}\right) \;. 
\end{equation}
By  Appendix~\ref{sec:Appendix},  for each character $\chi$ of
$\CA(N)$ we have 
\begin{equation}\label{uniform bound}
|E(\chi)|\leq 2 \sqrt{N} \;.
\end{equation}
%In our case, $D=2$. 

%Let $G=\langle A \rangle\subset\CA(N)$ be the group generated by $A$,
%whose size is $|G|=\ord(A,N)$. 
%Let $\chi$ be the multiplicative character of $G$ given by 
%$$ %\chi(g) = e(\frac{\tau}{\ord(A,N)}) ,\quad g=A^\tau \mod N\;.$$
%It is the restriction of 

%Define a multiplicative character $\chi_1$ of $\CA(N)$  as follows: 
Let $r= \frac{|\CA(N)|}{\ord(A,N)}$. 
%$r=[\CA(N):G] = \frac{|\CA(N)|}{\ord(A,N)}$.
Choose a generator $A_0$ of $\CA(N)$ such that
$A=A_0^r\mod N$. Define a character $\chi_1$ of $\CA(N)$ by setting
$\chi_1(A_0) = e\left(\frac 1{|\CA(N)|}\right)$.  
%Then $\chi_1(A) = e(\frac r{|\CA(N)|}) = e(\frac 1{\ord(A,N)}) = \chi(A)$.
Then 
$$
\chi_1(A^\tau) = e\left(\frac{\tau}{\ord(A,N)}\right) \;.
$$

%We define an ``incomplete'' version by 
%$$
%E_A(\chi) = \sum_{g\in G} \chi(g) e(\frac{F(g)}{N})
%$$

We may  write the indicator function of the subgroup of $\CA(N)$
generated by $A$ as
$$
{\mathbf 1}_{A}(y) = \frac 1{r} \sum_{\substack{\theta\in
\widehat{\CA(N)}\\ \theta(A)=1}} \theta(y)
$$ 
where the sum runs over all $r$ characters of $\CA(N)$ which are trivial on
$A$. Then we may rewrite $E_A(j)$ 
in terms of the complete sums \eqref{def complete sum} as
\begin{equation}\label{E(j) in terms of E(chi)} 
E_A(j) =  \frac 1{r}\sum_{\substack{\theta\in
\widehat{\CA(N)}\\ \theta(A)=1}} E(\chi_1^j\theta) \;.
\end{equation}
Now using the estimate \eqref{uniform bound} gives
$|E_A(j)|\leq 2\sqrt{N}$. 
\end{proof} 
%Remark: When $A$ generates all of the group $\CA(N)$, that is when
%$\ord(A,N)=N\pm 1$, only the trivial character $\chi=\mathbf 1$ 
%appears in \eqref{E(j) in terms of E(chi)} and in that case, the bound
%\eqref{uniform bound} is elementary, in fact $|E(\mathbf 1)| \leq 2$
%(see Appendix~\ref{sec:Appendix}).  

\vskip .5cm
\noindent{\bf Proof of Proposition~\ref{prop:S is small}} 
Let $k,k'\in \Z^2$, $Q(k)\neq Q(k')$.  Recall the definition 
\eqref{S general case} 
\begin{equation*}
S(k,k') = \sum_{t\neq 0\mod \ord(A,N)}  \Gamma(t) 
e\left(\frac{\overline{2}(q(k;A^t)-q(k';A^t))}{N}\right) \;.
\end{equation*}
%and  $$ F(g) = q(k;g)-q(k';g) $$
Expanding 
\begin{equation}\label{expand Gamma}
\Gamma(\tau) = \sum_{j\mod \ord(A,N)} \gamma(j) e\left(\frac {j\tau}{\ord(A,N)}\right) 
\end{equation} 
we get 
$$
S(k,k') = \sum_{j\mod \ord(A,N)} \gamma(j) E_A(j)
$$
where $E_A(j)$ is given in \eqref{defn E(j)}. 

According to Lemma~\ref{lem: bound E(j)},  if $Q(k)\neq Q(k')$ then  
$$
|S(k,k')|\leq 2\sqrt{N} \sum_{j\mod \ord(A,N)} |\gamma(j)|
$$
and so it remains to show that $\sum_j |\gamma(j)| = O(1)$. 

We first note that $\gamma(j) \geq 0$ so we may ignore the absolute
value signs: Indeed, from the definition \eqref{defn Gamma},
\eqref{expand Gamma}   we see that 
$$ \gamma(j) = \frac{L^2}{\ord(A,N)} \int_0^1
h_L(\theta)h_L\left( \theta+\frac j{\ord(A,N)}\right) d\theta
$$
which is non-negative since  $h_L\geq 0$. 

Thus we have 
$$\sum_j |\gamma(j)| = \sum_j \gamma(j) = \Gamma(0)$$
and by definition, 
$$
\Gamma(0) = \sum_{m\in \Z-\{0\}} \^h\left( \frac{\ord(A,N)}{L}m \right)^2
$$
which is bounded since $L<2\ord(A,N)$.
% and $\^h(x)\ll |x|^{-2}$. ????
\qed

\appendix
\section{An estimate for a character sum}\label{sec:Appendix}
In this appendix we give a proof for the bound 
$|E(\chi)|\leq2\sqrt{N}$ stated in  \eqref{uniform bound} 
for the character sum $E(\chi)$ defined in \eqref{def complete sum}.  
This bound is not new; 
as we explain below, in the ``split'' case it follows immediately from
Weil's bound \cite{Weil48}. In the ``inert'' case, the sum appears in
the work of Gurevich and Hadani \cite{GH} who discovered that the
matrix coefficients of $T_N(k)$ in the Hecke basis can be written as 
$$
\frac 1{\CA(N)} \sum_{B\in \CA(N)} \chi(B) 
\tr \left\{ T_N(k) \UN(B) \right\}
$$
and hence by \eqref{formula for q}, the matrix elements can be
expressed in terms of the sum $E(\chi)$. Gurevich and Hadani invoke
the full force of Deligne's Weil II paper \cite{Deligne Weil 2} to
give the bound  \eqref{uniform bound}. 
However, to make the paper more self contained, and perhaps also of
independent interest, we will give another proof that only requires
Weil's original methods \cite{Weil48}, together with some class field
theory.  Following Li \cite{Li92, Li}, 
we express the exponential sum in terms a certain id\`ele class
character sum (over degree one places), and then  derive the
bound from the Riemann Hypothesis for curves.
% This argument, originating in Li's work \cite{Li92} 
(The same argument was used in \cite{KurlbergAHP} in a similar
context.)

\subsection{$E(\chi)$ as a character sum} 
Using \eqref{formula for q} we can write $E(\chi)$ as follows: 
In the split case, where the matrix $A$ is diagonalizable over $\Z/N\Z$ then
$$
E(\chi) = \sum_{0,1\neq x\in \Z/N\Z} \chi(x)
\psi \left(\frac{1+x}{1-x} \right)
$$
where $\chi$ is a multiplicative character and $\psi$ a nontrivial 
additive character of $\Z/N\Z$. If $\chi\equiv \mathbf 1$ is trivial
then $E(\mathbf 1)=-\psi(1)-\psi(-1)$ so $|E(\mathbf 1)|\leq 2$. 
For $\chi\neq \mathbf 1$, the  bound $|E(\chi)|\leq 2\sqrt{N}$ follows
from Weil's 1948 result \cite{Weil48} 
(cf. \cite[Chapter 6, Theorem 3]{Li}). 

In the inert case, %we can write the sum $E(\chi)$ as follows: Let
let $\mathbb F$ be a quadratic extension of $\Z/N\Z$, $H\subset
\mathbb F^\times$ the group of elements of norm one, 
$\psi$ a nontrivial additive
character of $\Z/N\Z$, and $\chi$ a multiplicative character of
$H$. Let $\lambda_A\in H$, $\lambda_A\neq \pm 1$. Then by
\eqref{formula for q} 
\begin{equation*}
E(\chi) = \sum_{1\neq x\in H} \chi(x)
\psi\left(\frac1{\lambda_A-\lambda_A^{-1}} \frac{1+x}{1-x}\right) \;.
\end{equation*}

If the multiplicative character $\chi\equiv {\mathbf 1}$ is trivial, then
$E(\mathbf 1) = 0$, since 
$x\mapsto \frac 1{\lambda_A-\lambda_A^{-1}}\frac{1+x}{1-x}$ is a bijection of
$H\backslash \{1\}$ with $\Z/N\Z$.  From now on assume $\chi\neq
\mathbf 1$ is nontrivial. 

Take a quadratic non-residue $D\mod N$, and let 
$\sqrt{D}$ be a root of $X^2-D$ in $\mathbb F$. 
We may  write each element $1\neq x\in H$ uniquely as 
$$
x=\frac{t-\sqrt{D}}{t+\sqrt{D}}
$$
where $t\in \Z/N\Z$. In particular we
have $\lambda_A=\frac{t_0-\sqrt{D}}{t_0+\sqrt{D}}$ with $t_0\neq 0$
since $\lambda_A\neq \pm 1$. 
Then for $x\neq 1$, 
$$
\frac{1}{\lambda_A-\lambda_A^{-1}} \frac{1+x}{1-x} = 
\frac{D-t_0^2}{4t_0D }\cdot t
$$
and so 
$$
E(\chi) = \sum_{t\in \Z/N\Z}
\chi\left(\frac{t-\sqrt{D}}{t+\sqrt{D}}\right) 
\psi\left(\frac{D-t_0^2}{4t_0D }\cdot t \right) \;.
$$

Arguing as in \cite[Chapter 6]{Li}, we will  construct  id\`ele class
characters $\tilde \nu$, $\tilde \psi$ of the function field
$\Z/N\Z(X)$, of finite order,  satisfying: 

i) The conductors of $\tilde\psi$ and $\tilde \nu$ are 
$$ \cond(\tilde \psi) = 2\infty,\qquad \cond(\tilde \nu) = (w)$$
where $w$ is the degree two place of $\Z/N\Z(X)$  corresponding to the 
irreducible polynomial $X^2-D$. In particular the product $\tilde\psi
\tilde\nu$ is unramified at all finite degree one places $v\neq
\infty$. 
%field extension $\mathbb F$ of $\Z/N\Z$. 

ii)  Their values at a uniformizer $\pi_v$ for 
the degree-one place $v$ corresponding to the polynomial $X+t$ are 
$$
\tilde \psi(\pi_v) = \psi\left(\frac{D-t_0^2}{4t_0D }\cdot t \right) , \qquad 
\tilde \nu(\pi_v) = \chi\left(\frac{t-\sqrt{D}}{t+\sqrt{D}}\right) \;. 
$$
%(where $\pi_v$ is a uniformizer). 

Thus we can write $E(\chi)$ as a sum over  degree one 
places $v\neq \infty$  of $\Z/N\Z[X]$: 
$$
E(\chi) = \sum_{\substack{\deg(v)=1\\ v\neq \infty}}
(\tilde \nu \tilde \psi)(\pi_v)  \;.
$$

Class field theory and the Riemann Hypothesis for curves over a
function field give  (see \cite[Corollary 3 of Chapter 6]{Li})
$$
|E(\chi)|\leq (\deg\cond(\tilde\psi\tilde \nu) -2)\sqrt{N}\;. 
$$
Since the conductor of the product  $\tilde \nu \tilde \psi$ is
$2\infty +w$, which has degree $4$, we get 
$$
\left| E(\chi) \right| \leq 2\sqrt{N}
$$
as claimed. 

\subsection{Construction of id\`ele class characters}

We describe the construction of id\`ele class characters of the
function field  $K = \Z/N\Z(X)$. See \cite{Weil, Li} for background. 

Given a place $v$ of $K$, let $K_v$ denote the
completion of $K$ with respect to the topology induced by $v$, and
let $U_v = \{ \alpha \in K_v : |\alpha|_v = 1\}$ be the $v$-adic units
of $K_v$.  Let $\mathcal P_v=\{\alpha\in K_v: |\alpha|_v<1\}$ be the
maximal ideal, and denote by $\pi_v$ a uniformizer.  
In particular for the infinite place we may take $\pi_\infty = X^{-1}$. 

Let $I_K$ be  id\`ele group of $K$.  $I_K$ admits a product decomposition 
$I_K = K^\times \cdot \left( \prod_{v\neq
\infty} U_v \times K_\infty^\times \right)$,
with $\left(\prod_{v\neq\infty} U_v \times K_\infty^\times \right)
\cap K^\times = (\Z/N\Z)^\times$. 
The id\`ele class group is 
$$I_K/K^\times \simeq \left( K_\infty^\times \times \prod_{v\neq \infty} U_v
\right) / (\Z/N\Z)^\times \;.
$$  

\subsubsection{Constructing $\tilde \psi$} 
Given a nontrivial additive character $\psi_0$ of $\Z/N\Z$, we
will define an id\`ele class character $\tilde{\psi}$ %on $I_K/K^\times$
of finite order such that for the degree one place $v$ 
corresponding to the polynomial $X+t$ we have: 
$$
\tilde{\psi}( \pi_v) = \psi_0(t) \;.
$$ 
%Let $(0)$ be the place corresponding to the polynomial $X$. 
We first define $\tilde{\psi}$ on $U_\infty$ by setting
$$
\tilde{\psi}\left(a+bX^{-1}+\sum_{n\geq 2} c_n X^{-n}\right) 
%= \tilde{\psi}(a(1+(b/a)X)) 
= \psi_0(-b/a)
$$
so that we get a character of 
$U_\infty/\left((\Z/N\Z)^\times (1+\mathcal P_\infty^2)\right)$. 
Since $K_\infty^\times$ equals $\langle X^n\rangle_{n\in \Z} \times
U_\infty$, we may extend $\tilde\psi$ to a character of
$K_\infty^\times$ by declaring $\tilde \psi_\infty(X)=1$. 
Extend $\tilde{\psi}$ to $\prod_{v\neq \infty} U_v \times
K_\infty^\times$ by letting $\tilde{\psi}$ be trivial on 
$\prod_{v \neq \infty}U_v$.   
Since $\left(\prod_{v\neq\infty} U_v \times K_\infty^\times\right)
\cap K^\times = (\Z/N\Z)^\times$ and $\tilde{\psi}$ is trivial on
$(\Z/N\Z)^\times$, $\tilde{\psi}$ can be regarded as a character of
the id\`ele class group $I_K/K^\times$. 
%$I_K/K^\times \simeq \left( K_\infty^\times \times \prod_{v\neq \infty} U_v
%\right) / (\Z/N\Z)^\times $.  
%(note that the canonical map from  $\prod_{v\neq \infty} U_v \times
%K_\infty^*$ to $I_K/K^\times$ is  onto.)  

The conductor of $\tilde{\psi}$ is $2\cdot \infty$, since 
$\tilde{\psi}|_{U_v}$ is trivial for $v \neq \infty$ and
$\tilde{\psi}|_{1+p_\infty}$ is non-trivial.

Finally, the value of $\tilde{\psi}$ at the uniformizer $\pi_v$ for  a
degree one place $v$ corresponding to the polynomial $X+t$ equals  
\begin{multline*}
\tilde{\psi}_v(\pi_v) =
\tilde{\psi}_\infty\left(\frac 1{X+t}\right) = 
\tilde{\psi}_\infty(X\cdot (1+tX^{-1}))^{-1} \\
=\tilde{\psi}_\infty(1+tX^{-1})^{-1}  = \psi_0(-t)^{-1} = 
\psi_0(t) \;.
\end{multline*}

\subsubsection{Constructing $\tilde \nu$} 
Given a multiplicative character $\chi$ of the group $H$ of norm one
elements of $\mathbb F$, we define an
id\`ele class character $\tilde \nu$ of finite order so that 
$$
\tilde{\nu}( \pi_v) = \chi\left(\frac{t-\sqrt{D}}{t+\sqrt{D}}\right)
$$ 
if $v$ is
the degree one place corresponding to the polynomial $X+t$ (see
\cite[Chapter 6, proof of Theorem 6]{Li}): 
Denote by $w$ the degree two place corresponding to the irreducible
polynomial $X^2-D$.  
Then $U_w/(1+\mathcal P_w) \simeq \mathbb F^\times$ via the
map induced by $X\mapsto \sqrt{D}$. 
If $\sigma$ is the Galois involution of $\mathbb F$, then the map 
$x\mapsto \sigma(x)/x$ gives an
isomorphism  of $\mathbb F^\times/ (\Z/N\Z)^\times$ to the 
group $H\subset \mathbb F^\times$ of norm-one elements. 
Thus we get a homomorphism 
$$
\Phi:U_w\to U_w/\left( (\Z/N\Z)^\times\cdot
(1+\mathcal P_w)\right)\simeq H \;.
$$ 

%Given a multiplicative character $\chi$ of $H$, 
We  define a character $\tilde \nu_w$ of $U_w$ by 
$$
\tilde\nu_w(u):=\chi^{-1}(\Phi(u))\;, 
$$ 
which is trivial on $(\Z/N\Z)^\times
\left(1+\mathcal P_w \right)$. 
Extend it to a character $\tilde\nu$ of $\prod_{v\neq \infty} U_v
\times K_\infty^\times$ by having $\tilde\nu_v$ trivial if $v\neq w$. 
Since $\tilde \nu$ is trivial on 
$(\Z/N\Z)^\times =\left(\prod_{v\neq\infty} U_v \times
K_\infty^\times\right)\cap K^\times $,  $\tilde \nu$ gives a 
character of the id\`ele class group $I_K/K^\times$. 
%$I_K/K^\times \simeq \left( K_\infty^\times \times \prod_{v\neq \infty} U_v
%\right) / (\Z/N\Z)^\times $.  

If $\chi$ is nontrivial, then the conductor of $\tilde\nu$ is $w$. 
By construction, for a degree one place $v$ corresponding
to the polynomial $X+t$,  we have  
$$
\tilde\nu(\pi_v)= \tilde\nu_w\left(\frac 1{X+t}\right) = \chi(\Phi(X+t)) =
\chi\left(\frac{\sigma(\sqrt{D}+t)}{\sqrt{D}+t}\right)  =
\chi\left(\frac{t-\sqrt{D}}{t+\sqrt{D}}\right)  \;.
$$

%\bibliographystyle{amsplain}
%\bibliography{mybib,mypapers}

\end{document}